\documentclass[11pt,twoside]{article}
\usepackage{a4wide}
\usepackage{amssymb}
\usepackage{amsmath}
\newtheorem{theorem}{Theorem}

\newtheorem{corollary}{Corollary}
\newtheorem{lemma}{Lemma}

\newcommand{\lp}{\lambda_\psi}

\newcommand{\De}{\Delta}

\newcommand{\La}{\Lambda}

\renewcommand{\Bbb}[1]{\mathbb{#1}}
\newcommand{\N}{{\Bbb N}} % natural numbers
\newcommand{\R}{{\Bbb R}} % real numbers
\newcommand{\Rp}{{\Bbb R}^{+}} % positive real numbers
\newcommand{\Z}{{\Bbb Z}} % integer numbers

\newcommand{\cH}{{\cal H}}

\newcommand{\cR}{{\cal R}}

\newcommand{\cZ}{{\cal Z}}
\newcommand{\sabc}{\sigma_{a,b}(c)}
\newcommand{\sab}{\sigma_{a,b}}
\newcommand{\sabcd}{\sigma_{a',b'}(c')}
\newcommand{\sabd}{\sigma_{a',b'}}
\newcommand{\Zp}{\mathbb{Z}_{\geq0}}
\newcommand{\smz}{\setminus\{\mathbf{0}\}}
\newcommand{\hab}{h_{a,b}}

\newcommand{\ve}{\varepsilon}

\newcommand{\diam}{\operatorname{diam}}
\newcommand{\dist}{\operatorname{dist}}

\newcommand{\vv}[1]{{\mathbf{#1}}}

\newcommand{\ra}{R_{\alpha}}

\usepackage{verbatim}

\begin{document}

\title{Diophantine approximation with perfect squares and \\ the
solvability of an inhomogeneous wave equation}

\author {Victor Beresnevich\footnote{EPSRC Advanced Research Fellow, EP/C54076X/1}
\\ {\small\sc York} \and Maurice Dodson\footnote{Research supported by EPSRC
R90727/01} \\ {\small\sc York} \and Simon
Kristensen\footnote{William Gordon Seggie Brown Fellow}
\\ {\small\sc Edinburgh } \and Jason Levesley
\\ {\small\sc York } }

\maketitle

%\abstract{}

\bigskip

%%%\noindent{\small 2000 {\it Mathematics Subject Classification}\/:
%%%Primary ??; Secondary 11J83, 11J13, 11K60}\bigskip


\section{Introduction}

Diophantine criteria occur naturally in the theory of partial
differential equations through the notorious problem of small
denominators. An extensive treatment of such problems in the
theory of PDEs can be found, e.g., in \cite{Ptashnik84}. In this
paper, we are interested in a Diophantine problem related to an
inhomogeneous wave equation in $n$ spatial and one temporal
dimension with periodic boundary conditions. In brief, to ensure
the convergence of a formal solution to the equation certain
conditions on the periods should be satisfied. These conditions
normally leave a small set of exceptional periods for which the
convergence of the series is problematic, though the solution
might exist. It is therefore of interest to measure the `size' of
the exceptional set of periods. Regarding the wave equation we
will discuss the problem in more details and derive the
associated Diophantine problem in \S\ref{sec:problem}.

An analogous problem for the wave equation in one spatial
dimension is considered in \cite{Novak}. Even further, a more
general class of one dimensional PDEs is studied by Gramchev and
Yoshino in \cite{MR1355946}. However, their methods does not seem
to work in higher dimensions. In \cite{MR2004f:81068}, the
corresponding problem in two spatial dimensions is resolved for
the Schr\"odinger equation, for which the corresponding
Diophantine problem is partly linear, and it is settled by making
use of a result of Rynne \cite{MR93a:11066}.

The rest of the paper is structured as follows. The results of the
paper are stated in \S~\ref{sec:statements}. In
\S\S~\ref{sec:first_theorem}--\ref{sec:second_theorem} we prove
the results for the case when $n=2$ and in
\S~\ref{sec:outline_general} we outline how the proofs can be
adapted to obtain the $n$-dimensional versions.

Throughout we will use the Vinogradov notation: Given two real
valued functions $f$ and $g$, write $f \ll g$ if there is a
constant $c
> 0$ such that $f \le cg$. If $f \ll g$ and $g \ll f$, write
$f \asymp g$.

\section{The solubility of the wave equation and a related Diophantine
problem} \label{sec:problem}

Let $n \in \mathbb{N}$, $\alpha_i > 0$ for $i=1, \dots, n$,
$\beta>0$ and $f: \mathbb{R}^{n+1} \rightarrow \mathbb{R}$ be
periodic in all variables with period $\alpha_i$ in the $i$'th
variable and period $\beta$ in the $n+1$'st. We denote the $n$
first variables by $x_1, \dots, x_n$ and the $n+1$'st by $t$.
Suppose furthermore that $f$ is a smooth function of any of the
variables $x_i, t$, \emph{i.e.}, $f$ has continuous partial
derivatives of all orders. We will consider the partial
differential equation given by
\begin{equation}
\label{eq:1} \dfrac{\partial^2 u(\mathbf{x},t)}{\partial t^2} -
\Delta u(\mathbf{x},t) = f(\mathbf{x},t), \quad
\mathbf{x}=(x_1,\ldots,x_n) \in \mathbb{R}^n, t \in \mathbf{R},
\end{equation}
under the additional condition that the solution $u$ is smooth and
periodic with the same periods. Here $\Delta$ denotes the usual
Laplacian, \emph{i.e.},
\begin{equation*}
\Delta u(\mathbf{x},t) = \sum_{i=1}^n \dfrac{\partial^2
u(\mathbf{x},t)}{\partial x_i^2}.
\end{equation*}

The periodicity and smoothness conditions on $f$ are well-known to
be equivalent to the condition that $f$ has an expansion into a
Fourier series
\begin{displaymath}
f(\mathbf{x},t) = \sum_{(\mathbf{a},b) \in \mathbb{Z}^{n+1}}
f_{\mathbf{a},b} \exp\left(2 \pi i \left[\sum_{i=1}^n
\dfrac{a_i}{\alpha_i} x_i + \dfrac{b}{\beta} t \right] \right),
\end{displaymath}
where $\vv a=(a_1,\ldots,a_n)$, such that the coefficients
$f_{\mathbf{a},b}$ decay faster than the reciprocal of any
polynomial in $a_1, \dots, a_n ,b$ as $\max\{\vert a_1\vert,
\dots, \vert a_n \vert, \vert b \vert\}$ tends to infinity.

Suppose for the moment that \eqref{eq:1} has a solution $u$
satisfying the periodicity and smoothness conditions. Clearly, $u$
must also have the following Fourier expansion
\begin{displaymath}
u(\mathbf{x},t) = \sum_{(\mathbf{a},b) \in \mathbb{Z}^{n+1}}
u_{\mathbf{a},b} \exp\left(2 \pi i \left[\sum_{i=1}^n
\dfrac{a_i}{\alpha_i} x_i + \dfrac{b}{\beta} t \right] \right),
\end{displaymath}
Inserting this into \eqref{eq:1} and identifying coefficients, we
obtain
\begin{equation}
\label{eq:2} u_{\mathbf{a},b} = \dfrac{\beta^2}{4 \pi^2}
\dfrac{f_{\mathbf{a},b}}{\sum_{i=1}^n a_i^2
\frac{\beta^2}{\alpha_i^2} - b^2}.
\end{equation}
Now, since $\alpha_1, \dots, \alpha_n, \beta$ are fixed, and since
$f_{\mathbf{a},b}$ decays faster than the reciprocal of any
polynomial, for $u$ to be smooth it suffices to verify that
\begin{displaymath}
\left\vert \sum_{i=1}^n a_i^2 \frac{\beta^2}{\alpha_i^2} - b^2
\right\vert \geq C \max\{\vert a_1 \vert, \dots, \vert a_n
\vert\}^{-w},
\end{displaymath}
for some $C>0$, $w>1$ for all $(\mathbf{a},b) \in
\mathbb{Z}^{n+1}$ with $\vv a\not=\vv 0$. It is easy to see that
this condition can only fail if for any $w>1$ the inequality
\begin{equation}
\label{eq:3} \left\vert \sum_{i=1}^n a_i^2
\frac{\beta^2}{\alpha_i^2} - b^2 \right\vert < \max\{\vert a_1
\vert, \dots, \vert a_n \vert\}^{-w}
\end{equation}
holds for infinitely many $(\mathbf{a},b) \in \mathbb{Z}^{n+1}$
with $\vv a\not=\vv 0$.

Note that the condition given in \eqref{eq:3} is sufficient for
the solubility of \eqref{eq:1}, but not necessary. The Diophantine
problem considered in this paper is a natural generalisation of
the one of equation \eqref{eq:3}.

\section{Statement of results}
\label{sec:statements}

Throughout $\Zp$ will denote the set of non-negative integer
numbers and $|A|$ the Lebesgue measure of a set $A$. Given an
$n$-tuple $\mathbf{a}\in\Zp^2$, define the \emph{height}\/
$h_{\mathbf{a}}$ of $\mathbf{a}$ by setting
$h_\mathbf{a}:=\max(\,|\,a_1\,|\,\dots,|\,a_n\,|\,)$, that is
$h_{\mathbf{a}}$ is the highest coefficient of $\vv a$ in absolute
value.

Let $\psi:\Rp\to\Rp$ be a function such that $\psi(h)\to0$ as
$h\to\infty$ and define the set $W_n(\psi)$ to be
\begin{alignat*}{2}
W_n(\psi):=\{ \mathbf{x}\in[0,1]^n : &\ |\,
\mathbf{a}^2{}\cdot{}\mathbf{x} - b^2
% a^2x+b^2y-c^2
| <\psi(h_{\mathbf{a}}),
\\[1ex]
& \text{ holds for infinitely many }(\mathbf{a},b)\in\Zp^{n+1}\
\},
\end{alignat*}
where $\mathbf{a}^2 := (a_1^2,\dots,a_n^2)$.

The following statements constitute the main results of this
paper.
\begin{theorem}\label{t1}
Let $\psi:\Rp\to\Rp$ be monotonic. Then
\begin{equation*}
| W_n(\psi)|=
\begin{cases}
0,&
\sum_{h=1}^\infty h^{n-2}\psi(h)<\infty\,,\\[2ex]
1,&\sum_{h=1}^\infty h^{n-2}\psi(h)=\infty\,.
\end{cases}
\end{equation*}
\end{theorem}

\begin{theorem}\label{t2}
Let $\psi:\Rp\to\Rp$ be a monotonic. Given any positive $s< n$,
the $s$-dimensional Hausdorff measure of $W_n(\psi)$ satisfies the
relation
\begin{equation*}
\cH^s(W_n(\psi))=
\begin{cases}
0,& \sum_{h=1}^\infty \psi(h)^{s-(n-1)}h^{3n-2-2s}<\infty\,,\\[2ex]
\infty,&\sum_{h=1}^\infty \psi(h)^{s-(n-1)}h^{3n-2-2s}=\infty\,.
\end{cases}
\end{equation*}
\end{theorem}

\begin{corollary}\label{cor1}
Let $\psi:\Rp\to\Rp$ be a monotonic function such that
$\lim_{h\to\infty}\psi(h)=0$. Define $\lambda_\psi$, the lower
order of $1/\psi(2^r)$ at infinity, by setting
\begin{equation*}
\lambda_\psi=\liminf_{r\to\infty}\frac{-\log{\psi(2^r)}}{r\log 2}\
.
\end{equation*}
Note that $\lambda_\psi$ is always non-negative, but can be
infinity. If $n-1\le \lambda_\psi<\infty$ then
\begin{equation*}
\dim W_n(\psi)= (n-1)+\dfrac{n+1}{2+\lambda_\psi}\,.
\end{equation*}
In particular, if $\psi(r)=r^{-v}$ for some $v>n-1$ then
\begin{equation*}
\dim W_n(r\mapsto r^{-v})= (n-1)+\dfrac{n+1}{2+v}\,.
\end{equation*}
\end{corollary}

In terms of the wave equation, we may derive the following
corollary:
\begin{corollary}\label{cor2}
Let $\alpha_1, \dots, \alpha_n, \beta > 0$ and consider the
partial differential equation \eqref{eq:1}. Let $\delta_i =
\beta^2/\alpha_i^2$ for $i = 1, \dots, n$. If $f$ is smooth and
periodic in $x_1, \dots, x_n, t$ with periods $\alpha_1, \dots,
\alpha_n, \beta$ respectively, then \eqref{eq:1} is soluble with
$u$ smooth and periodic with the same periods whenever
$(\delta_1,\dots, \delta_n)$ does not belong to
\begin{equation*}
\bigcap_{v>1} W_n(r \mapsto r^{-v}),
\end{equation*}
a null set of Hausdorff dimension $n-1$.
\end{corollary}

\section{Proof of Theorem~\ref{t1}}
\label{sec:first_theorem}

We first prove the result for the case $n=2$ as the argument is
easiest to follow in this dimension.
%Then we give the outline
%of the proof for the general case when $n\geq3$ which is essentially
%the same.

\subsection{The case of convergence}

For every triple $(a,b,c)\in\Zp^3$ define the sets
$$
\sabc:=\{(x,y)\in[0,1]^2:|\,a^2x+b^2y-c^2|<\psi(\hab)\}
$$
$$
\sab:=\bigcup_{c\in\Z}\sabc.
$$
Without loss of generality we can assume that $a+b>0$. It is easy
to verify that
$$
| \sabc|\ll \frac{\psi(\hab)}{\hab^2}\,.
$$
Given a pair $(a,b)\in\Zp^2\setminus\{\vv0\}$,
$\sabc\not=\emptyset$ implies that $c\ll \hab$. It follows that
$$
| \sab|\ll \sum_{c\in\Zp\ :\
\sabc\not=\emptyset}\frac{\psi(\hab)}{\hab^2}\ll
\frac{\psi(\hab)}{\hab}\,.
$$
Now assume that $\sum_{h=1}^\infty\psi(h)<\infty$. Then,
\begin{equation}\label{e:001}
\sum_{h=1}^\infty \ \sum_{\substack{ (a,b)\in\Zp^2\smz;\\ \hab=h}}
| \sab|
\ll \sum_{h=1}^\infty \ \sum_{\substack{(a,b)\in\Zp^2\smz;\\
\hab=h}} \frac{\psi(h)}{h}\ll \sum_{h=1}^\infty\psi(h)<\infty.
\end{equation}
As the set $W_2(\psi)$ is exactly the set of points $(x,y)$ in the
unit square that fall into infinitely many sets $\sab$, we can
apply the Borel-Cantelli Lemma to (\ref{e:001}) to conclude that
the set $W_2(\psi)$ has zero Lebesgue measure.

\subsection{The case of divergence: Auxiliary Lemmas}

It should be noted that the main difficulty in proving
Theorem~\ref{t1} is in the case of divergence, to be considered in
sections~\ref{intersections} and \ref{divergence}. The line of
investigation of this case will rely on the following standard
auxiliary measure theoretic statements.

\begin{lemma}\label{lem1}
Let $\Omega$ be an open subset of $\R^n$ and let $|A|$ be the
Lebesgue measure of $A$. Let $E$ be a Borel subset of\/ $\R^n$.
Assume that there are constants $r_0,c>0$ such that for any ball
$B$ of radius $r(B)<r_0$ in $\Omega$  we have
$$|E\cap B|\ge c\ |B| \ . $$ Then $E$ has full measure in
$\Omega$, \textit{i.e.} $|\Omega\setminus E|=0$.
\end{lemma}

\begin{comment}
\begin{lemma}\label{lem1}
Let $\Omega$ be an open subset of $\R^n$ and let $\mu$ be the
Lebesgue measure in $\R^n$. Let $E$ be a Borel subset of\/ $\R^n$.
Assume that there are constants $r_0,c>0$ such that for any ball
$B$ of radius $r(B)
$$\mu(E\cap B)\ge c\ \mu(B) \ . $$ Then $E$ has full measure in
$\Omega$, \textit{i.e.} $\mu(\Omega\setminus E)=0$.
\end{lemma}
\end{comment}

\begin{lemma}\label{lem2}
Let $(\Omega,A,\mu) $ be a probability space and $E_n$ be a
sequence of $\mu$-measurable sets such that $\sum_{n=1}^\infty
\mu(E_n)=\infty $. Then $$ \mu( \limsup_{n \to \infty} E_n ) \;
\geq \; \limsup_{Q \to \infty} \frac{ \left( \sum_{s=1}^{Q}
\mu(E_s) \right)^2 }{ \sum_{s, t = 1}^{Q} \mu(E_s \cap E_t ) } \ \
\ . $$
\end{lemma}

In our particular problem we will take $E_n$ to be a subsequence
of the sequence of sets $\sab$. More precisely, we will estimate
pairwise intersections of $\sab$ restricted to a fixed ball $B$ on
average. The corresponding limsup set will be contained in
$W_2(\psi)\cap B$. On applying Lemma~\ref{lem2}, we will arrive at
a lower bound of the form $|W_2(\psi)\cap B|\ge c|B|$ for some
positive absolute constant. Lemma~\ref{lem1} will complete the
proof.

Further, to avoid painful and unnecessary calculation we will
restrict $B$ to be a ball lying inside $\Omega=[\ve,1]^2$ for some
arbitrarily small $\ve>0$. The corresponding probability measure
$\mu$ will be taken to be the normalized Lebesgue measure in
$\Omega$.

\subsection{Estimates for the measure of $\sab\cap B$ and their pairwise
intersections}\label{intersections}

Fix an arbitrary positive number $\ve<1$ and set
$\Omega:=[\ve,1]^2$. Take any ball $B$ in $\R^2$ lying in
$\Omega$.

\subsubsection{Restrictions on $c$}
\label{restrC}

Assume that $\sabc\cap B\not=\emptyset$. Then there is a point
$(x,y)\in B\subset[\ve,1]^2$ satisfying
$|a^2x+b^2y-c^2|<\psi(\hab)$. If $\hab$ is sufficiently large then
$\psi(\hab)<\ve$. Therefore, $ c^2<\ve+a^2x+b^2y\le 1+2\hab^2. $
Hence,
$$
| c|<2 \hab\,.
$$
On the other hand,
$$
c^2>a^2x+b^2y-\psi(h)>\ve(a^2+b^2)-\ve\ge\ve(h^2-1).
$$
Therefore,
$$
| c|>\ve\hab/2
$$
if $\hab$ is sufficiently large. Therefore, for all
$(a,b)\in\Zp^2$ with sufficiently large $\hab$ and all positive
$c$ with $\sabc\cap B\not=\emptyset$ we have
\begin{equation}\label{e:002}
\frac{\ve}{2}\,\hab<|c|<2\,\hab\ .
\end{equation}

\subsubsection{The amount of different $c$}

Define the line $R_{a,b,c}:=\{(x,y)\in\R^2:a^2x+b^2y-c^2=0\}$. It
is readily verified that $\sabc\cap B\not=\emptyset$ is equivalent
to $R_{a,b,c}\cap B\not=\emptyset$, except possibly for 2
`extremal' cases when $\sabc\cap B\not=\emptyset$ but the
corresponding lines do not hit the ball $B$ but lie sufficiently
close to $B$.

To evaluate the number of different $c$ such that
$\sabc\not=\emptyset$ we will estimate the number of lines
$R_{a,b,c}$ that hit the ball $B$ and then add 2 to the upper
estimate.

Let $x_0,y_0$ be the center of $B$ and $r$ be the radius of $B$.
Any point $(x,y)$ in $B$ can be written as
\begin{equation}\label{e:003}
x=x_0+\theta r\cos\phi,\quad y=y_0+\theta r\sin\phi,\qquad
0\le\theta<1,\ 0\le\phi<2\pi\,.
\end{equation}

Clearly, $R_{a,b,c}\cap B\not=\emptyset$ if and only if there is a
choice of $(x,y)$ subject to (\ref{e:003}) such that
$$
a^2(x_0+\theta r\cos\phi)+b^2(y_0+\theta r\sin\phi)-c^2=0.
$$
In such a case we have that
$$
c^2=a^2x_0+b^2y_0+\theta r(a^2\cos\phi+b^2\sin\phi)=
$$
$$
a^2x_0+b^2y_0+\theta
r\sqrt{a^4+b^4}\left(\frac{a^2}{\sqrt{a^4+b^4}}\cos\phi+\frac{b^2}{\sqrt{a^4+b^4}}\sin\phi\right)=
$$
$$
a^2x_0+b^2y_0+\theta
r\sqrt{a^4+b^4}\left(\sin\phi_0\cos\phi+\cos\phi_0\sin\phi\right)=
$$
$$
a^2x_0+b^2y_0+\theta r\sqrt{a^4+b^4}\,\sin(\phi+\phi_0),
$$
where $\phi_0=\arcsin\frac{a^2}{\sqrt{a^4+b^4}}$. Therefore, $c^2$
varies in the interval
\begin{equation}\label{e:004}
[a^2x_0+b^2y_0-r\sqrt{a^4+b^4},\ a^2x_0+b^2y_0+r\sqrt{a^4+b^4}].
\end{equation}
Moreover, on taking $\phi:=\pm\pi/2-\phi_0$,
$\sin(\phi+\phi_0)=\sin(\pm\pi/2)=\pm1$ we see that any perfect
squares in this interval does contribute to a line $R_{a,b,c}$
which hits the ball $B$. Clearly $c^2$ lies in (\ref{e:004}) if
and only if $c$ is in the interval
\begin{equation}\label{e:004sub1}
\left[\ \sqrt{a^2x_0+b^2y_0-r\sqrt{a^4+b^4}},\
\sqrt{a^2x_0+b^2y_0+r\sqrt{a^4+b^4}}\ \right]\ .
\end{equation}
The length of interval (\ref{e:004sub1}) is
$$
\xi_{a,b,B}=\sqrt{a^2x_0+b^2y_0+r\sqrt{a^4+b^4}}-\sqrt{a^2x_0+b^2y_0-r\sqrt{a^4+b^4}}=
$$
$$
\frac{2r\sqrt{a^4+b^4}}{\sqrt{a^2x_0+b^2y_0+r\sqrt{a^4+b^4}}+\sqrt{a^2x_0+b^2y_0-r\sqrt{a^4+b^4}}}.
$$
Taking into account that $\ve\le x_0,y_0\le1$ and $r<1$, it
follows that
$$
\frac12\ r\,\hab\le \xi_{a,b,B}\le \frac{8}{\ve}\ r\,\hab\,.
$$
Now, the number of possible values for $c$ lies between
$\xi_{a,b,B}$ and $\xi_{a,b,B}+3$ and is therefore $\asymp r\,
\hab$.

\subsubsection{The measure of $\sab\cap B$}
\label{sec:estsigcapB}

Given a $c$, it is easily verified that $|\sabc\cap B|\le 4
r\psi(\hab)/\sqrt{a^4+b^4}\le 4 r\psi(\hab)/\hab^2$, where $r$ is
the radius of $B$.

The number of possible values of $c$ such that $\sabc\cap
B\not=\emptyset$ is bounded above by $\xi_{a,b,B}+3\le
\frac{10}\ve\ r\, \hab$ if $\hab$ is sufficiently large.
Therefore,
$$
| \sab\cap B|\le 4 r\psi(\hab)/\hab^2\times \frac{10}\ve\ r\, \hab
= c_2|B|\ \frac{\psi(\hab)}{\hab},
$$
where $c_2=\frac{40}{\ve\pi}$ and $\hab$ is sufficiently large.

Let $\frac12B$ be the ball centred at the same point as $B$ of
radius $r/2$. Then it is an elementary geometric task to compute
that $|\sabc\cap B|\ge r\psi(\hab)/\hab^2$ whenever $\sabc\cap
\frac12B\not=\emptyset$ and $\hab$ is sufficiently large.

The number of possible values of $c$ such that $\sabc\cap
\frac12B\not=\emptyset$ is bounded below by $\xi_{a,b,\frac12B}\ge
\frac14 \ r\,\hab$. Therefore,
$$
| \sab\cap B|\ge r\psi(\hab)/\hab^2 \times \frac14 \
r\,\hab=c_1|B|\frac{\psi(\hab)}{\hab},
$$
where $c_1=\frac{1}{4\pi}$.

The upshot of the above is that
\begin{equation}\label{e:005}
c_1|B|\frac{\psi(\hab)}{\hab} \le |\sab\cap B|\le
c_2|B|\frac{\psi(\hab)}{\hab}
\end{equation}
for all sufficiently large $\hab$, where $c_1,c_2$ are absolute
positive constants.

\subsubsection{Additional conditions on $(a,b)$}
\label{sect:add_conds_on_ab}

Throughout the remainder of the proof of Theorem~\ref{t1} we will
assume that the following conditions on $(a,b)$ hold:
\begin{equation}\label{e:006}
\gcd(a,b)=1,
\end{equation}
where $\gcd$ means the greatest common divisor, and
\begin{equation}\label{e:007}
1/2\le a/b\le 2.
\end{equation}
The above conditions sift elements of the sequence of sets $\sab$
which prevent us from having sufficiently good estimates for the
measures of pairwise intersections of these sets. On the other
hand, the remaining `thinned out' part of the sequence $\sab$ is
still rich enough to ensure that the sum
\begin{equation}\label{e:008}
\sum|\sab|
\end{equation}
diverges over this restricted sequence. Such a condition as that
of Equation (\ref{e:008}) is necessary to apply Lemma~\ref{lem2}.
Indeed, to verify that (\ref{e:008}) diverges over $(a,b)\in\Zp^2$
satisfying (\ref{e:006}) and (\ref{e:007}) define $N_k$ to be the
number of $(a,b)$ satisfying (\ref{e:006}) and (\ref{e:007}) with
$2^k\le \hab<2^{k+1}$. Then in view of symmetry of the set of
$(a,b)$ of interest we get
$$
N_k=2\sum_{2^k\le a<2^{k+1}}\ \sum_{\stackrel{\scriptstyle
\rule[-1ex]{0ex}{2ex} b<a}{\text{(\ref{e:006}) and (\ref{e:007})
are satisfied}}}1= 2\sum_{2^k\le a<2^{k+1}}\
\Big(\varphi(a)-\varphi([a/2])\Big),
$$
\begin{comment}
$$ N_k=2\sum_{2^k\le a<2^{k+1}}\ \sum_{\stackrel{\scriptstyle
\rule[-1ex]{0ex}{2ex} b satisfied}}}1= 2\sum_{2^k\le a<2^{k+1}}\
\Big(\varphi(a)-\varphi([a/2])\Big),
$$
\end{comment}
where $\varphi$ is the Euler function. It is well known that
$$
\sum_{1\le q\le Q}\varphi(q)=\frac{3}{\pi^2}Q^2+O(Q\log Q).
$$
Then
$$
2\sum_{2^k\le a<2^{k+1}}\
\varphi(a)=\frac6{\pi^2}((2^{k+1})^2-(2^k)^2)+O(k2^k)=\frac{18}{\pi^2}
\ 2^{2k}+O(k2^k)
$$
and
$$
2\sum_{2^k\le a<2^{k+1}}\ \varphi([a/2])=4\sum_{2^{k-1}\le
x<2^{k}}\
\varphi(x)=\frac{12}{\pi^2}((2^{k})^2-(2^{k-1})^2)+O(k2^k)=\frac9{\pi^2}
\ 2^{2k}+O(k2^k).
$$
It follows that
$$
N_k=\frac9{\pi^2} \ 2^{2k}+O(k2^k).
$$
Now the estimated sum is
$$
\sum_{\stackrel{\scriptstyle \rule[-1ex]{0ex}{2ex}
(a,b)\in\Zp^2}{\text{(\ref{e:006}) and (\ref{e:007}) are
satisfied}}}|\sab\cap B|=\sum_{k=0}^\infty\sum_{2^k\le h<2^{k+1}}
\sum_{\stackrel{\scriptstyle
\rule[-1ex]{0ex}{2ex}(a,b)\in\Zp^2\,:\,
\hab=h}{\text{(\ref{e:006}) and (\ref{e:007}) are satisfied}}} |
\sab|\gg
$$
$$
\gg |B|\sum_{k=0}^\infty\sum_{2^k\le h<2^{k+1}}
\sum_{\stackrel{\scriptstyle
\rule[-1ex]{0ex}{2ex}(a,b)\in\Zp^2\,:\,
\hab=h}{\text{(\ref{e:006}) and (\ref{e:007}) are satisfied}}}
\frac{\psi(2^{k+1})}{2^k}\asymp
$$
$$
\asymp|B|\sum_{k=0}^\infty2^k\psi(2^k)\asymp|B|
\sum_{h=1}^\infty\psi(h)=\infty.
$$

Finally, note that the limsup set for the `thinned out' sequence
$\sab$ is contained in the limsup set for the complete sequence
$\sab$, which is $W_2(\psi)$. Therefore, it will be sufficient to
prove that the thinned out limsup set is of full Lebesgue measure
in order to ensure that $W_2(\psi)$ is also of full measure.

An immediate consequence of condition (\ref{e:006}) is that for
any two pairs $(a,b)$ and $(a',b')$ satisfying (\ref{e:006}) the
assumption $(a,b)\not=(a',b')$ implies that $(a,b)$ and $(a',b')$
are not collinear. Moreover, $(a^2,b^2)$ and $(a'^2,b'^2)$ are not
collinear. Therefore we can assume that the (smaller) angle
between $(a^2,b^2)$ and $(a'^2,b'^2)$, which will be denoted by
$\alpha=\alpha(a,b,a',b')$, is not zero. The analysis of the
measures of intersections $\sab\cap\sabd\cap B$ will rely on the
behaviour of this angle and is given in the following sections.

\subsubsection{The measure of intersections in the case of a big angle}
\label{sec:big_angle}

We will assume that $(a,b)\not=(a',b')$. Within this subsection we
set $h=\hab$ and $h'=h_{a',b'}$. For simplicity we will assume
that $h\ge h'$. Now
\begin{equation}\label{e:009}
\sab\cap\sabd\cap B=\bigcup_{c'\in\Zp}\sab\cap\sabcd\cap B
\end{equation}
For a fixed $c'$ the set $\sabcd\cap B$ is covered with a strip of
length $2r$ (recall that $r$ is the radius of $B$) and width
$\psi(h')/h'^2$. This strip is a piece of the
$\psi(h')/h'^2$-neighbourhood of the line
\begin{equation}\label{e:010}
a'^2x+b'^2y-c'^2=0.
\end{equation}
To estimate the measure in (\ref{e:009}) we first estimate the
measure of the intersection of $\sab$ with such a strip.

The angle $\alpha=\alpha(a,b,a',b')$ introduced in the previous
section is the (smaller) angle between the line defined in
(\ref{e:010}) and the family of parallel lines
\begin{equation}\label{e:011}
a^2x+b^2y-c^2=0, \text{ where $c\in\Zp$.}
\end{equation}
Using (\ref{e:002}) it is readily verified that the distance
between two consecutive lines in the family (\ref{e:011}) is
$\asymp h^{-1}$.

Now if $A$ and $B$ are two consecutive points on the line
(\ref{e:010}) obtained as a result of its intersection with two
consecutive lines in (\ref{e:011}), say $\mathcal{L}_1$ and
$\mathcal{L}_2$, it is easy to calculate that the distance between
$A$ and $B$ is the distance between $\mathcal{L}_1$ and
$\mathcal{L}_2$ divided by $\sin\alpha$, that is $\asymp
\frac{1}{h\sin\alpha}$. Since the piece of the line (\ref{e:010})
of interest is of length at most $2r$, there are at most
$$
\ll rh\sin\alpha+1
$$
non-empty intersections $\sabc\cap\sabcd\cap B$ when $c$ runs over
all integers.

As the set $\sabc\cap\sabcd$ is a parallelepiped with area $\ll
\frac{\psi(h)}{h^2}\,\frac{\psi(h')}{h'^2}\,\frac{1}{\sin\alpha}$,
the upshot of the above is that
$$
| \sab\cap\sabcd\cap B|\ll
\frac{\psi(h)}{h^2}\,\frac{\psi(h')}{h'^2}\,\frac{1}{\sin\alpha}\times(rh\sin\alpha+1).
$$
Further, since there are $\ll rh'$ values of $c'$ that need to be
considered, we have that
\begin{equation}\label{c}
\begin{array}[b]{rl}
| \sab\cap\sabd\cap B|\ll &\displaystyle
\frac{\psi(h)}{h^2}\,\frac{\psi(h')}{h'^2}\,\frac{1}{\sin\alpha}\times(rh\sin\alpha+1)rh'\asymp\\[3ex]
&\displaystyle\asymp|B|\,\frac{\psi(h)}{h}\,\frac{\psi(h')}{h'}\,\left(1+\frac1{rh\sin\alpha}\right).
\end{array}
\end{equation}

Assuming that $\frac1{rh\sin\alpha}\le 1$, or equivalently that
\begin{equation}\label{e:012}
\sin\alpha\ge \frac1{rh},
\end{equation}
gives
\begin{equation}\label{e:013}
| \sab\cap\sabd\cap B|\ll |
B|\,\frac{\psi(h)}{h}\,\frac{\psi(h')}{h'}\,.
\end{equation}
Finally, since there are $\asymp h$ integer vectors $(a,b)$ with
$\hab=h$ and $\asymp h'$ integer vectors $(a',b')$ with
$h_{a',b'}=h'$, summing the measures of intersections
$|\sab\cap\sabd\cap B|$ in the case under consideration results in
$$
\sum_{\stackrel{\rule[-1ex]{0ex}{3ex}\scriptstyle\hab\le H,\
h_{a',b'}\le H}{(a,b)\not=(a',b')\text{ and (\ref{e:012})
holds}}}|\sab\cap\sabd\cap B|\ll |
B|\left(\sum_{h=1}^H\psi(h)\right)^2.
$$

\subsubsection{The measure of intersections in the case of a small angle}
\label{sec:small_angle}

In this section we will deal with the case of
\begin{equation}\label{e:014}
\sin\alpha\le \frac1{rh}\,.
\end{equation}
Again we will assume that $(a,b)\not=(a',b')$ and given a matrix
$A$, $|A|$ will denote its determinant and $\|A\|$ the absolute
value of its determinant.

Since $\alpha$ is the angle between the vectors $(a^2,b^2)$ and
$(a'^2,b'^2)$ it follows that
\begin{equation}\label{e:015}
h^2h'^2\sin\alpha\asymp\sqrt{a^4+b^4}\sqrt{a'^4+b'^4}\
\sin\alpha=\left\|\begin{array}{cc}a^2&b^2\\
a'^2&b'^2\end{array}\right\|= \left\|\begin{array}{cc}a&b\\
a'&b'\end{array}\right\|\times \left\|\begin{array}{cc}a&-b\\
a'&b'\end{array}\right\|\ .
\end{equation}
If $\beta$ denotes the (smaller) angle between $(a,b)$ and
$(a,-b)$ then
$$
\sin\beta=\frac1{a^2+b^2}\ \left\|\begin{array}{cc}a&b\\
a&-b\end{array}\right\|=\frac{2|ab|}{a^2+b^2}\
\stackrel{(\ref{e:007})}{\ge}\ \frac12\,.
$$
Hence, $\beta\ge\pi/6$ and the angle between $(a',b')$ and at
least one of the vectors $(a,b)$ and $(a,-b)$ is at least
$\pi/12$. Without loss of generality we can assume that such an
angle is between $(a,-b)$ and $(a',b')$. Then
$$
\left\|\begin{array}{cc}a&-b\\
a'&b'\end{array}\right\|=\sqrt{a^2+b^2}\sqrt{a'^2+b'^2}\sin\pi/12\gg
h\, h'\,.
$$
It now follows from (\ref{e:015}) that
\begin{equation}\label{e:016}
1\le
\left\|\begin{array}{cc}a&b\\
a'&b'\end{array}\right\|\ll h\, h'\,\sin\alpha\le\frac{h'}{r}\,.
\end{equation}
This means that for every fixed $a',b',a$ there are at most
$\ll\frac1r$ possible values for $b$. Indeed, $|ab'-a'b|\ll
h'r^{-1}$, that is $|b-ab'/a'|\ll h'r^{-1}/a'\ll r^{-1}$.
Moreover, (\ref{e:016}) implies that
\begin{equation}\label{e:017}
\sin\alpha\gg\frac1{h\,h'}\,.
\end{equation}
To complete the analysis for this case we consider two specific
subcases.

\noindent{\bf Subcase (i)~--~moderately small angle.}~\\
Assume for the moment that
\begin{equation}\label{cc}
  \sin\alpha\ge\frac1{r^2\,h\,h'}\,.
\end{equation}
Using (\ref{c}), (\ref{e:014}) and (\ref{cc}) it follows that
$$
| \sab\cap\sabd\cap B|\ll |
B|\,\frac{\psi(h)}{h}\,\frac{\psi(h')}{h'}\,\frac1{rh\sin\alpha}\ll
| B|\,\frac{\psi(h)}{h}\,\psi(h')\ r\ .
$$
Now the sum of intersections for this subcase can be estimated as
follows,
$$
\sum_{\hab\le H,\ h_{a',b'}\le H} |\sab\cap\sabd\cap B|\ll
$$
$$
\sum_{h=1}^H\sum_{h'=1}^{h-1}\sum_{h_{a',b'}=h'}\sum_{\hab=h} |
\sab\cap\sabd\cap B|\ll
$$
$$
\sum_{h=1}^H\sum_{h'=1}^{h-1}\sum_{h_{a',b'}=h'}\sum_{\hab=h} |
B|\,\frac{\psi(h)}{h}\,\psi(h')\ r \ \ll
$$
\begin{equation}\label{v}
\sum_{h=1}^H\sum_{h'=1}^{h-1}h' |B|\,\frac{\psi(h)}{h}\,\psi(h')
\ll | B|\sum_{h=1}^H\sum_{h'=1}^{h-1} \psi(h)\,\psi(h') \ll |
B|\left(\sum_{h=1}^H \psi(h)\right)^2\,.
\end{equation}

\noindent{\bf Subcase (ii)~--~ultra small angle.}~\\
To complete the analysis of all possible values of $\alpha$ it
remains to consider the case when
$$
\sin\alpha<\frac1{r^2\,h\,h'}\,.
$$
Then
\begin{equation}\label{cv}
\left\|\begin{array}{cc}a&b\\
a'&b'\end{array}\right\|\ll h\, h'\,\sin\alpha\le\frac{1}{r^2}\,.
\end{equation}
and
\begin{equation}\label{cvc}
| \sab\cap\sabd\cap B|\ll |
B|\,\frac{\psi(h)}{h}\,\frac{\psi(h')}{h'}\,\frac1{rh\sin\alpha}\ll
| B|\,\frac{\psi(h)}{h}\,\psi(h')\frac1r.
\end{equation}
Now we estimate the number of quadruples $(a,b,a',b')$ satisfying
(\ref{e:006}), (\ref{e:007}), (\ref{cv}), $2^k\le\hab<2^{k+1}$ and
$2^l\le h_{a',b'}<2^{l+1}$. Given fixed $a$ and $b'$, (\ref{cv})
means that $a'$, $b$ can only be chosen to satisfy $|ab'-a'b|\ll
r^{-2}$. This means that there are $\ll r^{-2}$ possible values
for $t=a'b$. In turn, for a fixed $t$ there are at most $d(t)$
possible values for $a'$ and $b$, where $d(t)$ is the number of
divisors of $t$. It is well known that for any $\delta>0$ there is
a constant $c_{\delta}>0$ such that $d(t)\le c_\delta t^\delta$
for all $t$. Taking $\delta=1/4$ we get that the number of
possible quadruples $a,b,a',b'$ is $\ll (2^k\ 2^l)^{5/4}r^{-2}$.

Without loss of generality we assume that $\psi(h)\le h^{-1}$.
Then the sum of intersections for this subcase is estimate as
follows
$$
\sum_{\hab\le H,\ h_{a',b'}\le H} |\sab\cap\sabd\cap B|=
$$
$$
\sum_{k=1}^{[\log H]+1}\ \ \sum_{l=1}^{[\log H]+1}\sum_{2^k\le
\hab< 2^{k+1},\ 2^l\le h_{a',b'}< 2^{l+1}} |\sab\cap\sabd\cap
B|\ll
$$
$$
\sum_{k=1}^{[\log H]+1}\ \ \sum_{l=1}^{k} |
B|\,\frac{\psi(2^k)}{2^k}\,\psi(2^l)\frac1r \times (2^k\
2^l)^{5/4}r^{-2}\ll \frac1r \sum_{k=1}^{[\log H]+1}\
\sum_{l=1}^{k} 2^{k/4}\psi(2^k)2^{5l/4}\psi(2^l)\ \ll
$$
$$
\frac1r\sum_{k=1}^{[\log H]+1}\ \sum_{l=1}^{k}
2^{3k/4}\psi(2^k)2^{3l/4}\psi(2^l)\ \ll \frac1r\sum_{k=1}^\infty\
\sum_{l=1}^\infty 2^{-k/4}2^{-l/4} <\infty.
$$
We are now in a position to complete the proof of Theorem~\ref{t1}
for the divergence case.

\subsection{Completion of the proof of Theorem~\ref{t1}}\label{divergence}

The upshot of the above computations is the following estimates:
$$
S_1(H)=\sum_{(a,b)\in\cZ_H}|\sab\cap B|\gg |
B|\left(\sum_{h=1}^H\psi(h)\right)
$$
$$
S_2(H)=\sum_{(a,b)\in\cZ_H}\ \ \
\sum_{(a',b')\in\cZ_H}|\sab\cap\sabd\cap B|\ll |
B|\left(\sum_{h=1}^H\psi(h)\right)^2
$$
where $\cZ_H=\{(a,b)\in\Zp^2,\ \text{(\ref{e:006}) and
(\ref{e:007}) hold and }\ \hab\le H\}$. Therefore,
$$
\frac{S_1(H)^2}{S_2(H)}\gg |B|
$$
for all sufficiently large $H$. Since
$\limsup_{\hab\to\infty}\sab\cap B\subset W_2(\psi)\cap B$, by
Lemma~\ref{lem2}
$$
| W_2(\psi)\cap B|\ge |\limsup_{\hab\to\infty}\sab\cap B|\gg |B|.
$$
This holds for any ball $B$ in $\Omega$ with the implied constant
independent of $B$. Therefore, by Lemma~\ref{lem1}, $W_2(\psi)$
has full measure in $\Omega=(\ve,1)^2$. Since $\ve>0$ is
arbitrary, $W_2(\psi)$ has full measure in $[0,1]^2$. This
completes the proof of Theorem~\ref{t1}.

\section{Proof of Theorem~\ref{t2}}
\label{sec:second_theorem}

\subsection{Hausdorff measures and dimension}

In this section we give a very brief introduction to the theory of
Hausdorff measures and dimension. For further details consult
\cite{MAT}.

Let $s$ be a positive real number. The Hausdorff $s$--measure will
be denoted throughout by $\cH^{s}$ and is defined as follows.
Suppose $F$ is a non--empty subset of $\R^k$. Suppose that
$\rho>0$. A $\rho$-cover of $F$ is a countable collection $\{ B_i
\}$ of balls in $\R^k$ with radii $r_i \leq \rho $ for each $i$
such that
$$
F \subset \bigcup_{i} B_{i}.
$$
Define the function $\cH^s_\rho$ by
$$
\cH^s_\rho (F) := \inf \left\{ \sum_{i} r_i^s \right\}
$$
where the infimum is taken over all possible $\rho$-covers of $F$.
Then $\cH^s(F)$ of the set $F$ is defined by
$$
\mathcal{ H}^{s} (F) := \lim_{ \rho \to0} \cH^{s}_{\rho} (F) \; =
\; \sup_{\rho > 0 } \mathcal{ H}^{s}_{\rho} (F) \; .
$$
Let $F$ be an infinite set. The Hausdorff dimension of $F$ is the
(unique) number
$$
\dim F=\inf\{s>0:\cH^s(F)=0\}=\sup\{s>0:\cH^s(F)=+\infty\}.
$$
Note that $\cH^k$ is a multiple of the $k$-dimensional Lebesgue
measure in $\R^k$ when $k\in\N$.

\subsection{Proof of Theorem~\ref{t2}. The case of convergence}

The proof of convergence is straightforward.
%Throughout this section we set $h=h_{a,b}$.
Recall from above that $W_2(\psi)$ can be expressed as a limsup
set of the form
$$
W_2(\psi) =
\bigcap_{h=1}^\infty\bigcup_{\substack{(a,b)\in\Z^2,\\
h_{a,b}=h}}^\infty\bigcup_{c\in\Z} \left(\sigma_{a,b}(c)\right).
$$
Each $\sigma_{a,b}(c)$ can be covered by a family $C_{a,b}^c$ of
balls each of radius $\psi(h_{a,b})/h^2_{a,b}$ where
$$
\sharp{}C_{a,b}^c \ll \frac{h_{a,b}^2}{\psi(h_{a,b})}.
%\{c : \sigma_{a,b}(c)\cap\mathbb{I}^2 \neq \emptyset \} \ll
%h_{a,b}^2.
$$
By assumption $\psi(h)\to 0$ as $h\to\infty$. Therefore, given any
$N\in\N $, $\psi(h)/h^2\leq 1/N$ for sufficiently large $h$. It
follows that
\begin{align*}
\cH^s_{1/N}(W_2(\psi)) &\ll \sum_{\substack{(a,b)\in\Z^2,\\
h_{a,b}\geq{}N}} \left( \frac{\psi(h_{a,b})}{h_{a,b}^2} \right)^s
\frac{h_{a,b}^2}{\psi(h_{a,b})} h_{a,b} \ll \sum_{h\geq{}N} \left(
\frac{\psi(h)}{h^2}
\right)^s \psi(h)^{-1} h^2 h h \\
&= \sum_{h\geq{}N} \psi(h)^{s-1}h^{4-2s} \to 0 \text{ as }
N\to\infty.
\end{align*}
Therefore $\cH^s(W_2(\psi))=0$, as required.

\subsection{Proof of Theorem~\ref{t2}. The case of divergence}

To prove the divergence case of Theorem~\ref{t2} we appeal to a
recent result of Beresnevich \& Velani \cite{BersVelSlice} in
which a mass transference principle for linear forms based on a
technique called `slicing' is established. The result allows one
to transfer statements about the Lebesgue measure of general
limsup sets occurring in Diophantine approximation to ones
involving Hausdorff measure.

The ideas outlined below are specialised to suit the particular
Diophantine approximation problems posed in this paper and are
therefore simplified versions of those given in
\cite{BersVelSlice}. The general framework of \cite{BersVelSlice}
is far richer and allows one to address Diophantine problems
involving systems of linear forms, inhomogeneous approximation and
general measure functions in one consuming package.

%Let $k, m \geq 1 $ and $ l\ge0$ be integers with $k=m+l$ and
Let $\cR=(\ra )_{\alpha \in J}$ be a family of lines in $\R^2$
indexed by an infinite countable set $J$. For every $\alpha\in J$
and $\delta\geq 0$ define the $\delta$--neighborhood $
\Delta(R_\alpha,\delta)$ of $R_\alpha$ by
$$
\Delta(R_\alpha,\delta) := \{\vv x \in \R^2: \dist(\vv x,R_\alpha)
< \delta\} \ .
$$
%where $\dist(\vv x, R_\alpha) = ????$
%$$ \Delta(R_\alpha,\delta) := \{\vv x \in \R^k: \dist(\vv
%x,R_\alpha) < \delta\} \ . $$ Thus $\Delta(R_\alpha,\delta)$ is
%simply the $\delta$--neighborhood of the $l$--dimensional plane
%$R_\alpha$. Note that by definition, $ \Delta(R_\alpha,\delta) =
%\emptyset $ if $\delta =0$.
Next, let
$$\Upsilon : J \to \R^+ : \alpha\mapsto
\Upsilon(\alpha):=\Upsilon_\alpha$$ be a non-negative, real valued
function on $J$.
%In order to avoid pathological situations within
%our framework, we
Further, assume that for every $\epsilon
> 0$ the set $\{\alpha\in J:\Upsilon_\alpha>\epsilon \}$ is finite.
This condition implies that $\Upsilon_\alpha \to 0 $ as $\alpha$
runs through $J$. Now define the following `$\limsup$' set,
$$ \La(\Upsilon)=\{\vv x\in\R^2:\vv
x\in\Delta(R_\alpha,\Upsilon_\alpha)\ \mbox{for\ infinitely\ many\
}\alpha\in J\} \ . $$
%Note that in view of the conditions imposed on
%$k,l$ and $m$ we have that $l < k$. Thus the dimension of the
%`approximating' planes $ R_\alpha$ is strictly less than that of the
%ambient space $\R^k$. The situation when $l =k $ is of little
%interest and has therefore been naturally omitted.

\begin{theorem}\label{BVslice}
Let $\cR$ and $\Upsilon$ as above be given. Let $V$ be a line in\/
$\R^2$ and

~ \qquad \qquad $(i)$\quad \ $V \ \cap \ R_\alpha \ \neq \
\emptyset $ \quad for all $ \ \alpha\in J \ $,

~ \qquad \qquad $(ii)$\quad $\sup_{\alpha\in J}\diam( \,
V\cap\De(R_\alpha,1) \, ) \ < \ \infty \ $ .

\noindent Let $f$ and $g : r \to g(r):= r^{-1} \, f(r)$ be
dimension functions such that $r^{-2}f(r)$ is monotonic and let
$\Omega $ be a ball in $\R^2$. Suppose for any ball $B$ in
$\Omega$
$$
\cH^2 \big(\, B \cap \La \big(g(\Upsilon) %^{\frac1m}
\big) \, \big) \, = \, \cH^2(B)
$$ Then
$$ \cH^f \big( \, B \cap \Lambda(\Upsilon) \, \big) \, = \, \cH^f(B) \ . $$
\end{theorem}

Now, let $f : r\to r^s$. As $1 < s < 2$ it follows that
$r^{-2}f(r)$ is monotonic and $f$ and $g$, defined as above, are
both dimension functions. Further, let $\Omega$ to be the unit
square $[0,1)^2$, $J:=\{ (a,b,c)\in\Z^3: h_{a,b}=|a|\}$,
% ?? SHOW HOW WE KNOW $h_{a,b}=|a|$ IS ENOUGH
$$
R_{(a,b,c)} = \{ (x,y)\in\R^2: a^2{}x+b^2{}y = c^2 \}
$$
and $\Upsilon_{(a,b,c)} := \psi(h_{a,b})/h_{a,b}^2$. Define
$S_2(\psi)$ to be
$$
S_2(\psi) := \Lambda(\Upsilon)\cap[0,1)^2.
$$
Note that $S_2(\psi)\subset{}W_2(\psi)$ and $|S_2(\psi)|=1$
whenever $|W_2(\psi)|=1$. To complete the proof of Theorem
\ref{t2}, it is sufficient to prove the divergence case for
$S_2(\psi)$. With this in mind, let $V:=\{(x_1,x_2)\in\R^2: x_2=0
\}$. It is straightforward to verify that conditions (i) and (ii)
of Theorem \ref{BVslice} hold in this case. From the divergence
case of Theorem \ref{t1}, it follows that $\cH^2(S_2(\psi)) = 1 =
\cH^2([0,1)^2) $. Therefore,
$\cH^s(S_2(\psi))=\cH^s(\mathbb{I}^2)=\infty$ and Theorem \ref{t2}
is proved.

\subsection{Proof of Corollary~\ref{cor1}}

By the definition of the lower order for any $\delta>0$ the
inequality
$\lambda_\psi+\delta\ge\frac{\log\frac1{\psi(2^r)}}{\log 2^r}$ for
infinitely many $r$. It follows that
\begin{equation}\label{e:021}
\psi(2^r)\ge (2^r)^{-\lambda_{\psi}-\delta}\text{ \ \ for
infinitely many }r\,.
\end{equation}
Take $s=1+\frac{3}{2+\lp+\delta}-\delta$. Then
$$
\psi(2^r)^{s-1}(2^r)^{5-2s}\ge
(2^r)^{-(\lambda_{\psi}+\delta)(s-1)+5-2s}=
(2^r)^{-(\lambda_{\psi}+2+\delta)(s-1)+3}=
(2^r)^{\delta(\lp+2+\delta)}>1
$$
for infinitely many $r$. Therefore,
$$
\sum_{r=1}^\infty\psi(2^r)^{s-1}(2^r)^{5-2s}=\infty\,.
$$
Since $\psi$ is monotonic, using a simple `condensation' argument
it is easy to verify that
$$
\sum_{h=1}^\infty\psi(h)^{s-1}h^{4-2s}=\infty\,.
$$
Hence, by Theorem~\ref{t2},
$$
\cH^s(W_2(\psi))=\infty\qquad\text{ and }\qquad\dim W_2(\psi)\ge
s=1+\frac{3}{2+\lp+\delta}-\delta.
$$
Since $\delta>0$ is arbitrary, we have $\dim W_2(\psi)\ge
1+\frac{3}{2+\lp}$.

Again, by the definition of the lower order, for any $\delta>0$
the inequality
$\lambda_\psi-\delta\le\frac{\log\frac1{\psi(2^r)}}{\log 2^r}$
holds for all sufficiently large $r$. It follows that
\begin{equation}\label{e:022}
\psi(2^r)\le (2^r)^{-\lambda_{\psi}+\delta}\text{ \ \ for all
sufficiently large }r\,.
\end{equation}
Take $s=1+\frac{3}{2+\lp-\delta}+\delta$. Then
$$
\psi(2^r)^{s-1}(2^r)^{5-2s}\le
(2^r)^{-(\lambda_{\psi}-\delta)(s-1)+5-2s}=
(2^r)^{-\delta(\lp+2-\delta)}
$$
for infinitely many $r$. Therefore,
$$
\sum_{r=1}^\infty\psi(2^r)^{s-1}(2^r)^{5-2s}<
\sum_{r=1}^\infty(2^r)^{-\delta(\lp+2-\delta)}<\infty\,.
$$
Since $\psi$ is monotonic, using the `condensation' argument it is
easy to verify that
$$
\sum_{h=1}^\infty\psi(h)^{s-1}h^{4-2s}<\infty\,.
$$
Hence
$$
\cH^s(W_2(\psi))<\infty\qquad\text{ and }\qquad\dim W_2(\psi)\le
s=1+\frac{3}{2+\lp-\delta}+\delta.
$$
Since $\delta>0$ is arbitrary, we have $\dim W_2(\psi)\le
1+\frac{3}{2+\lp}$. Therefore, we have the equality $\dim
W_2(\psi) = 1+\frac{3}{2+\lp}$.

\subsection{Proof of Corollary~\ref{cor2}}

%??

The proof that Equation (\ref{e:001}) has a solution in
$\mathbb{H}^{m+2}(\alpha,\beta,\gamma)$ whenever
$(\delta_1,\delta_2) \notin W_2(r \mapsto r^{-2})$, which is a set
of dimension $7/4$, is an immediate consequence of
Corollary~\ref{cor1}.

Assume now that $f$ is required to be smooth. As $W_2(r \mapsto
r^{-\tau^\prime}) \subset W_2(r \mapsto r^{-\tau})$ for
$\tau^\prime
> \tau$. It follows by continuity of $\dim(\cdot)$ that
$$
\mathrm{dim}\left( \bigcap_{v>1} W_2(r \mapsto r^{-v}) \right) =
\lim_{v\to\infty} \dim\left(W_2(r \mapsto r^{-v})\right) =
\lim_{v\to\infty} \left(1+\frac{3}{2+v}\right)=1.
$$
This establishes Corollary~\ref{cor2}

\section{Outline of the General case $n\geq3$}
\label{sec:outline_general}

The convergence case of Theorem~\ref{t1} for $n\geq3$ is almost
immediate. For every $(n+1)$-tuple $(\mathbf{a},b)\in\Zp^{n+1}$,
let
$$
\sigma_\mathbf{a}(b) := \{x\in[0,1]^n : |\mathbf{a}^2\cdot{}x -
b^2 | < \psi(h_\mathbf{a})\}
$$
and
$$
\sigma_\mathbf{a} := \bigcap_{b\in\Z} \sigma_\mathbf{a}(b)
$$
where $\mathbf{a}^2$ is the vector $(a_1^2,a_2^2,\dots,a_n^2)$. It
is easy to see that each set $\sigma_\mathbf{a}(b)$ is an
$n-1$-dimensional hyperplane with area $|\sigma_\mathbf{a}(b)| \ll
\psi(h_\mathbf{a})/h^2_\mathbf{a}$. Fix an $\mathbf{a}\in\Zp^n$,
$\sigma_{\mathbf{a}}\neq\emptyset$ implies that
$b\ll{}h_\mathbf{a}$. Note that the number of vectors $\mathbf{a}$
for which $h_\mathbf{a}=h$ is $\ll{}h^{n-1}$. Now
$$
\sum_{h=1}^\infty\sum_{\substack{ \mathbf{a}\in\Zp^n\smz;\\
h_\mathbf{a}=h}} \sum_{\substack{ b\in\Z:\\
\sigma_\mathbf{b}\neq\emptyset}}|\sigma_\mathbf{a}(b) | \ll
\sum_{h=1}^\infty{}h^{n-2}\psi(h) < \infty
$$
by assumption. It follows that $|W_n(\psi)|=0$ and we are done.

Assuming for a moment the validity of the divergence part of
Theorem~\ref{t1} when $n\geq3$. Establishing Theorem~\ref{t2} is
relatively straightforward.

In the convergence case we note that
$$
W_n(\psi) =
\bigcap_{h=1}^\infty\bigcup_{\substack{\mathbf{a}\in\Z^n,\\
h_\mathbf{a}=h}}^\infty\bigcup_{b\in\Z}
\left(\sigma_\mathbf{a}(b)\cap\mathbb{I}^n\right)
$$
and each $\sigma_\mathbf{a}(b)$ can be covered by a family
$C^b_\mathbf{a}$ of balls each of radius
$\psi(h_\mathbf{a})/h^2_\mathbf{a}$ such that
$$
\#{}C^b_\mathbf{a} \ll \left( h_\mathbf{a}^2 /\psi(h_\mathbf{a})
\right)^{n-1}.
$$
It is then a simple matter to amend the proof in the case when
$n=2$ for $n\geq3$ and deduce that $\mathcal{H}^s(W_n(\psi))=0$.

The divergence case of Theorem~\ref{t2} can be proved with only
minor modifications of the proof for the case when $n=2$. The main
changes to be made to the general framework of
Theorem~\ref{BVslice} are that $\mathcal{R}$ is now a countable
family of $(n-1)$-dimensional hyperplanes, $\mathbf{x}\in\R^n$,
$V$ is a linear subspace of $\R^n$, $f$ is a dimension function
such that $r^{-n}f(r)$ is monotonic and $g:r\to{}r^{-(n-1)}f(r)$
is a dimension function.

Now, let $f:r\to{}r^s$, $\Omega$ be the unit hypercube $[0,1)^n$,
$J:=\{ (\mathbf{a},b)\in\Z^{n+1}_{\ge0}:h_\mathbf{a}=|a_1|\}$,
$$
R_{(\mathbf{a},b)} := \{ \mathbf{x}\in\R^n: \mathbf{a}^2\cdot{}\vv
x = b^2 \}
$$
and $\Upsilon_{(\mathbf{a},b)} :=
\psi(h_\mathbf{a})/h_\mathbf{a}^2. $ The rest of the argument is
essentially the same as that given above with $2$ replaced by $n$
and $V := \{ \mathbf{x}\in\R^n: x_n=0 \}$.

It remains to establish the divergence part of Theorem~\ref{t1}
for the cases when $n\geq{}3$. As noted above, the family of lines
that we considered in \S~\ref{sec:first_theorem} have now been
replaced by $(n-1)$-dimensional hyperplanes, but the analysis
again hinges on the angle between the members of two non-collinear
families. It is relatively easy to see that the restrictions that
applied to $c$ in \S~{\ref{restrC}} must also apply to $b$ in the
above argument and further, that the number of such $b$ must also
be $\asymp{}rh_\mathbf{a}$. This follows from the fact that the
geometry in the $n$-dimensional case can be reduced to the same
problem as that of the $2$-dimensional case by projecting the ball
$B$ and the $(n-1)$-dimensional hyperplanes onto a $2$-dimensional
plane perpendicular to the family of hyperplanes defined by the
equations
$$
\mathbf{a}^2\cdot\mathbf{x} - b^2 = 0
$$
where $b\in\Z$. A simple geometric argument implies that
$|\sigma_\mathbf{a}(b)\cap{}B| \ll
r^{n-1}\displaystyle\frac{\psi(h_\mathbf{a})}{h_\mathbf{a}^2}$
where $r$ is the radius of $B$. As the number of possible $b$ such
that $\sigma_\mathbf{a}(b)\cap{}B\neq\emptyset$ is $\ll r
h_\mathbf{a}$ it follows that
$$
| \sigma_\mathbf{a}\cap{}B| \ll
r^n\frac{\psi(h_\mathbf{a})}{h_\mathbf{a}^2}{}h_\mathbf{a} \ll |
B|\frac{\psi(h_\mathbf{a})}{h_\mathbf{a}},
$$
and by an analogous argument to that in \S~\ref{sec:estsigcapB} it
can be shown that
$$
| \sigma_\mathbf{a}\cap{}B| \gg |
B|\frac{\psi(h_\mathbf{a})}{h_\mathbf{a}}
$$
where the constants implied by the $\ll$ and $\gg$ are absolute.
Recall that conditions~(\ref{e:006}) and ~(\ref{e:007}) were
imposed on $a$ and $b$ in the $2$-dimensional cases. For the
higher dimensional cases the corresponding conditions become
\begin{equation}
    \label{e:023}
    \mathrm{gcd}(a_1,a_2,\dots,a_n) = 1
\end{equation}
and
\begin{equation}
    \label{e:024}
    1/2 \leq a_1/a_2 \leq 2,
\end{equation}
with the same consequences as in \S~\ref{sect:add_conds_on_ab},
namely a sufficient quantity of vectors to maintain divergence of
our sum and non-collinearity of any two vectors
satisfying~(\ref{e:023}).

As in the $2$-dimensional case considered above, take any two
vectors $\mathbf{a}$ and $\mathbf{a}^\prime$ with
$\mathbf{a}\neq\mathbf{a}^\prime$, which must be linearly
independent by (\ref{e:023}). The upshot of linear independence is
that the angle between the normals to the two hyperplanes, and
therefore the hyperplanes themselves, is non-zero. Strictly
speaking there are two angles, but we shall take the smaller of
the two and call this $\alpha$. The result of
\S~\ref{sec:big_angle} also holds in this case. It is a simple
geometric argument to show that the volume of the parallelepiped
obtained by intersecting any two members of the two families is
now
$$\ll
r^{n-2}\frac{\psi(h_\mathbf{a})}{h_\mathbf{a}^2}\frac{\psi(h_{{\mathbf{a}}^\prime})}{h_{{\mathbf{a}}^\prime}^2}.$$
An analogous argument to that presented in \S~\ref{sec:big_angle}
with the restriction that $\sin\alpha\geq\frac{1}{rh}$ yields the
desired estimate for the sum of the measures of the intersections
subject to the above restriction on $\alpha$.

To complete the proof requires taking care of the cases when the
angle $\alpha$ becomes small. Recall that in the $2$-dimensional
case, \S~\ref{sec:small_angle},  this naturally split into two
cases; that of a moderately small angle and an ultra-small angle.
It was shown in the former case that the same estimate as that of
the big angle case could be deduced and in the latter, that the
sum of the intersections over the class of vectors with ultra
small angle was in fact convergent and could therefore be
neglected. It is precisely these conclusions that can be shown to
hold in the general case and the divergence part of
Theorem~\ref{t1} will follow in exactly the same manner as in the
$2$-dimensional case.

The analysis in \S~\ref{sec:small_angle} relied on a key
observation that the angle, $\alpha$, couldn't get too small. More
precisely that $\sin\alpha\gg
1/h_\mathbf{a}h_{\mathbf{a}^\prime}$. This was a consequence of
the assumption that $1/2\leq a_1/a_2 \leq 2$. To establish this
fact we used the standard result from elementary geometry that
$|\mathbf{a}\times\mathbf{b}| = |
\mathbf{a}||\mathbf{b}||\sin\beta|$ where $\beta$ is the angle
between $\mathbf{a}$ and $\mathbf{b}$. In higher dimensions the
cross product $\times$ is replaced by the wedge product $\wedge$
where
$$
    \mathbf{a}\wedge\mathbf{b} = \left\{ \ \left|\begin{array}{cc}a_i & a_j
    \\[1ex]
b_i & b_j \end{array}\right|\ :\  1\leq{}i<{}j\leq{}n \right\}.
$$
Note without any loss of generality we can assume that the first
two coordinates give the biggest determinant by reordering if
necessary and it is this observation, coupled with the assumption
that $1/2\leq a_1/a_2\leq{}2$ that allows us to conclude that
$\sin\alpha\gg 1/h_{\mathbf{a}}h_{\mathbf{a}^\prime}$. The
argument for the case when the angle is moderately small is
exactly the same as for the $2$-dimensional case. Leaving only the
case when
\begin{equation}
\label{eq:100}
\sin\alpha<\frac{1}{r^2h_{\mathbf{a}}h_{\mathbf{a}^\prime}}
\end{equation}
to take care of. As there is a free choice in all but the first
two components of either of the vectors $\mathbf{a}$ and
$\mathbf{a}^\prime$ the number of pairs of vectors that we need to
consider is
$h_\mathbf{a}^{n-2}h_{\mathbf{a}^\prime}^{n-2}\cdot\#\{(a_1,a_2,a_1^\prime,a_2^\prime)
\}$. Using the estimate we deduced in \S~\ref{sec:small_angle} it
follows that the sum we are estimating is convergent and can
therefore be neglected.

The final steps in proving the divergence part of Theorem~\ref{t1}
follow in exactly the same manner as that of the $2$-dimensional
case.

There are only minor modifications needed to the proofs of
Corollaries~\ref{cor1} and \ref{cor2} to establish them in the
general case and the details are left to the reader.

%$$
%\sum_{\stackrel{\rule[-1ex]{0ex}{3ex}\scriptstyle h_{\mathbf{a}}\le
%H,\ h_{\mathbf{a}^\prime}\le H}{\mathbf{a}\not=\mathbf{a}^\prime
%\text{ and } \sin\alpha\geq1/rh_\mathbf{a}}} |\s\cap\sabd\cap B|\ll
%|B|\left(\sum_{h=1}^H h^{n-2}\psi(h)\right)^2.
%$$

%---- plain!, alpha, unsrt, abbrv, siam!, amsalpha

\end{document}